\documentclass[preprint,12pt]{elsarticle}
\usepackage{amsmath}
\usepackage{amssymb}
\usepackage{amsfonts}
\usepackage{eucal}
\usepackage[usenames]{color}
\usepackage{graphicx}
\numberwithin{equation}{section}
\newtheorem{Theorem}{Theorem}[section]

\newtheorem{Definition}[Theorem]{Definition}

\newtheorem{Remark}[Theorem]{Remark}

\begin{document}
\title{The arctangent law for a certain random time related to a one-dimensional diffusion}
\author{Mario Abundo }
\address{{Dipartimento di Matematica, Universit\`a  ``Tor Vergata'', via della Ricerca Scientifica, I-00133 Rome,
Italy.\\
E-mail: \tt{abundo@mat.uniroma2.it}}

\begin{abstract}
\noindent  For a time-homogeneous, one-dimensional diffusion process $X(t),$
we investigate the distribution of the first instant, after  a given time $r,$ at which
$X(t)$  exceeds its maximum on the interval $[0,r],$
generalizing a result of Papanicolaou,  which is valid for
Brownian motion.
\end{abstract}

}
\date{}
\maketitle

\noindent {\bf Keywords:} One-dimensional diffusion, First-passage time, Maximum value on an interval \\
{\bf Mathematics Subject Classification:} 60J60, 60H05, 60H10.

\section{Introduction}
In this paper, we extend to a one-dimensional diffusion process $X(t)$
the result of (Papanicolaou, 2016) for  Brownian motion, concerning the arctangent law for a
certain random time.\par\noindent
Indeed, let be $X(t)$ a time-homogeneous, one-dimensional diffusion in the interval $I \subset \mathbb{R}$ which is the solution of the SDE:
\begin{equation} \label{eqdiffu}
dX(t) = \mu (X(t)) dt + \sigma (X(t)) d B_t , \ X(0) = \eta \in I,
\end{equation}
where $B_t $ (with $B_0=0)$ is standard Brownian motion (BM) and the drift and diffusion
coefficients satisfy the usual conditions (see e.g. (Ikeda and Watanabe, 1981)) for existence and uniqueness of
the solution of \eqref{eqdiffu}.  \par\noindent
For a fixed time $r >0,$ we consider the maximum $M_r := \max _{ 0\le t \le r } X(t) $ of the diffusion $X$ on the
interval $[0,r],$ and we denote by $S$ the following random time:
\begin{equation} \label{Sdefinition}
S= S (r) := \inf \{ t \ge r: X(t)\ge  M_r \} -r .
\end{equation}
Assuming that the initial state $\eta$ is random, our aim is to study the distribution function of $S,$
generalizing the result of (Papanicolaou, 2016), that refers to
the case when $X(t)= B^ \eta _t := \eta + B_t $ (i.e. BM starting from the random value  $ \eta,$ not necessarily zero),  and states that:
\begin{equation} \label{distrtauBM}
P \{ S _B (r) \le t \} = \frac 2 \pi \arctan \left ( \sqrt {\frac t r } \ \right ), \ t \ge 0,
\end{equation}
where $S_B(r) = \inf \{ t \ge r : B^ \eta _t \ge \max _ {s \in [0,r]} B^ \eta _s \} -r .$
By taking the derivative with respect to $t,$ one obtains the probability density of
$S _B (r) :$
\begin{equation} \label{densityofSB}
f_ {S _B (r) } (t)= \frac {\sqrt r } { \pi (r+t) \sqrt t} , \ t \ge 0 .
\end{equation}
Notice that the expectation, $E( S _B (r)),$ turns out to be
infinite. \par The knowledge of the distribution of $S$ is
relevant in various diffusion models used in applied sciences,
such as Mathematical Finance, Biology, Physics, Hydraulics, etc.,
whenever  the time evolution of the phenomenon under study is
described by a diffusion $X(t);$ in fact,  one is often interested
to find the first instant, after a given time $r,$ at which $X(t)$
exceeds the maximum  value attained in the time interval $[0,r],$
namely in times prior to $r.$ For instance, in the Economy
framework, if we let $r$ vary in $(0, + \infty),$ the process
$S(r),$ so obtained, is related to the drawdown  process, which
measures the fall in value of $X(t)$ from its running maxima, and
is frequently used as performance indicator in the fund management
industry (see e.g. (Dassios and Lim, 2017) and references
therein). Indeed, $S(r)$ can be expressed in terms of the time
elapsed since the last time the maximum is achieved, that was
studied in (Dassios and Lim, 2017).

\section{The result}
Let $w(x) \in C^2 (I) $ be the {\it scale function} associated to the diffusion $X(t)$ driven by the SDE \eqref{eqdiffu},
that is, the  solution of:
\begin{equation} \label{scaleeq}
\begin{cases}
L w (x) = 0 , \ x \in I  \\
  w(0) = 0, \ w '(0) = 1,
\end{cases}
\end{equation}
where $L$ is the infinitesimal generator of $X$ defined by:
\begin{equation} \label{generator}
L h = \frac 12 \sigma ^2 (x) \frac {d^2 h } {d x^2}  + \mu (x) \frac {d h } {d x}, \ \
h \in C^2 (I) .
\end{equation}
Actually, the scale function can be taken as any
function $\widetilde w = aw + b,$ with $a > 0$ and $b \in \mathbb{R}$ (see e.g. (Karlin and Taylor, 1975)); we chose
the initial conditions of \eqref{scaleeq}, for the sake of simplicity. \par\noindent
As easily seen,  if the integral $ \int _0 ^t \frac { 2 \mu (z)}
{\sigma ^2 (z)} \ dz $ converges, the problem \eqref{scaleeq}  has solution:
\begin{equation} \label{scalefunction}
w(x)= \int _ 0 ^x \exp \left ( - \int _0  ^t \frac { 2 \mu (z)}  {\sigma ^2 (z)} \ dz \right ) dt .
\end{equation}
If $\zeta (t) := w(X(t)),$ by It${\rm \hat o}$'s formula one obtains
\begin{equation}
 \zeta (t) = w ( \eta) + \int _0 ^t w '( w ^{-1} (\zeta (s)) ) \sigma ( w ^{-1} (\zeta  (s))) d B_s  \ ,
\end{equation}
that is, the process $\zeta (t) $ is a local
martingale, whose quadratic variation is
\begin{equation} \label{rho}
\rho (t) \doteq \langle \zeta  \rangle _t = \int _0 ^ t [ w' (X(s)) \sigma (X(s)) ] ^ 2 ds, \ t \ge 0 .
\end{equation}
The (random) function $\rho(t)$ is differentiable, increasing, and
$\rho (0)=0.$ If $\rho ( + \infty) = + \infty,$ it can be shown (see e.g. (Revuz and Yor, 1991)) that there exists a
BM $\widehat B$ such that $\zeta (t) = \widehat B ( \rho (t)) + w ( \eta );$
thus, since $w$ is invertible, the solution $X(t)$ to \eqref{eqdiffu} can be written in the form
\begin{equation} \label{representation}
X(t)= w^ {-1} ( \widehat B(\rho (t)) + w (\eta) ) .
\end{equation}
In this way, $X$ is obtained from BM by a space transformation and a random time-change (see e.g. the discussion in (Abundo, 2012)).

\begin{Definition} {\rm
We say that the diffusion $X(t)$  (with $X(0)= x )$ is {\it conjugated } to BM (see also (Abundo, 2012)),
if there exists an increasing differentiable function $v(x)$ with $v(0) = 0,$ such that
$X(t)= v^{-1} \left (B_t + v(x ) \right ),$ for any $t \ge 0.$
}
\end{Definition}

\begin{Remark} {\rm Diffusions conjugated to BM are special cases of \eqref{representation}, for
$\rho (t)=t, \ \widehat B_t = B_t,$ and $w=v$ (however, it is not required that  $v'(0)=1)$.
}
\end{Remark}

A class of diffusions  conjugated to BM is given by processes $X(t)$ which are solutions of SDEs such as:
\begin{equation} \label{conjdiffu}
 dX(t) = \frac 1 2 \sigma (X(t)) \sigma ' (X(t)) dt + \sigma (X(t)) dB_t, \  X(0)=x,
\end{equation}
with $\sigma (\cdot) \ge 0.$
Indeed, if the integral
$v(x) :=  \int  ^ x \frac 1  {\sigma (r) }dr $
is convergent,
by It${\rm \hat o}$'s  formula, one obtains $X(t)= v^ {-1} (B_t +v(x)).$ \par\noindent
Explicit examples of diffusions conjugated to BM (see also (Abundo, 2012)) are: \bigskip

\noindent
$\bullet$ the diffusion in $I = \rm I\!R$ driven by the SDE $dX(t)= \frac 1 3 X(t) ^ {1/3} dt + X(t) ^ {2/3} dB_t, \ X(0)= x ,$ which is conjugated to BM via the  the function $v(x)= 3 x^ {1/3},$ that is,
$X(t)= (x ^ {1/3} + \frac 1 3 B_t ) ^3 ;$ \par\noindent
$\bullet$  the diffusion in $I= [0, + \infty )$ driven by the SDE $dX(t)= \frac 1 4 dt + \sqrt { X(t) \vee 0} \ dB_t , \ X(0)= x \ge 0 $
(Feller process), which is conjugated to BM
via the  the function $v(x)= 2 \sqrt x,$ that is,
$X(t)= \frac 1 4 (B_t + 2 \sqrt x ) ^2 ;$ \par\noindent
$\bullet$   the diffusion in $I= [0, 1]$  driven by the SDE \par\noindent
$dX(t)= (\frac 1 4 - \frac 1 2 X(t))dt + \sqrt { X(t) (1- X(t)) \vee 0 } \ dB_t , \ X(0)= x \in [0,1] $ (Wright-Fisher like process), which  is conjugated to BM
via the  the function $v(x)= 2 \arcsin \sqrt x,$ that is,
$X(t)= \sin ^2 ( B_t/2 + \arcsin \sqrt x ) .$ \bigskip

\noindent
If we drop the requirement that $v(0)=0$ in Definition 2.1, then,
for $\sigma >0,$ the diffusion in $I= (0, + \infty )$  driven by the SDE  $dX(t)= \frac { \sigma ^2 } 2 X(t) dt + \sigma X(t) dB_t , \  X(0)=x >0 $ (a special case of geometric BM)
is conjugated to BM
via the function $v(x)= \frac {\ln x } \sigma,$ that is, $X(t)= \exp ( \sigma B_t + \ln x ).$
\bigskip

\noindent The class of processes $X$ given by  \eqref{representation} with $\rho (t)$ deterministic, includes, besides diffusions conjugated to BM, the integral of Gauss-Markov processes (see (Abundo, 2015), (Abundo, 2013)), e.g.
integrated BM  \ $X(t)= \int _0 ^t B_s ds , $ represented by $X(t) = \widehat B (\rho (t))$ with $\rho (t) = t^3/3 , $
and  integrated Ornstein-Uhlenbeck (OU) process $X(t) = \int _0 ^t Y(s) ds,$ where $Y(t)$ is OU process
(see  (Abundo, 2013) for the explicit representation of $X$ in the form \eqref{representation}).

\bigskip

The announced result is:

\begin{Theorem} \label{teorema}
Let $X(t)$ be the solution of the SDE \eqref{eqdiffu} and suppose that the scale function $w$ of $X,$ given by \eqref{scaleeq} exists; for fixed $r>0,$ let  $S(r)$ be the random time defined
by \eqref{Sdefinition}.
With the previous notations, suppose that the function $\rho$ satisfies the condition $\rho (+ \infty ) = + \infty;$ \par\noindent
(i) if $\rho (t)$ is deterministic, then the probability distribution of $S(r)$ is:
\begin{equation} \label{result}
P \{ S (r) \le t \} = \frac 2 \pi \arctan \left ( \sqrt {\frac {\rho (t+r) - \rho (r)} {\rho (r)} } \ \right ), \ t \ge 0,
\end{equation}
and its density is:
\begin{equation} \label{denresult}
f_{ S (r)} (t)= \frac { \rho '(t+r) \sqrt {\rho (r)} } {\rho (t+r) \sqrt {\rho (t+r) - \rho (r)} } .
\end{equation}
(ii) If $\rho$ is not deterministic, let us suppose that there exist two deterministic, continuous increasing
functions $\alpha (t)$ and $ \beta (t),$ with $\alpha (0) = \beta (0) =0,$ such that
\begin{equation} \label{boundsrho}
\alpha (t) \le \rho (t) \le \beta (t),  \ \forall t \ge 0.
\end{equation}
Then:
\begin{equation} \label{bound2}
P \{ S(r) \le t \} \le \frac 2 \pi \arctan \left ( \sqrt {\frac {\beta (t+r) - \alpha (r)} {\alpha (r)} } \ \right ), \ t \ge 0.
\end{equation}
Moreover, if there exists $\bar t>0$ such that $\alpha (r+ \bar t) \ge \beta (r) ,$ then:
\begin{equation} \label{bound1}
P \{ S(r) \le t \} \ge \frac 2 \pi \arctan \left ( \sqrt {\frac {\alpha (t+r) - \beta (r)} {\beta (r)} } \ \right ) ,
\  t >  \bar t  .
\end{equation}
\end{Theorem}
{\it Proof.} \ Under the hypothesis, the representation \eqref{representation} of $X$ holds. \par\noindent
(i) Suppose that $\rho (t)$ is
deterministic. Then:
$$ M_r = \max _ { t \in [0,r]} X(t) =  \max _ {t \in [0,r]} w ^{-1} \left ( \widehat B ( \rho (t)) + w(\eta) \right )
= w^{-1} \left ( \max _ {t \in [0, \rho (r)]} \widehat B_t + w(\eta) \right ),$$
and so:
$$ \tau _r: = \inf \Big \{ t \ge r: X(t)\ge  M_r \Big \}
= \inf \Big \{ t \ge r: \widehat B ( \rho (t))   \ge \max _ { u \in [0, \rho (r)] } \widehat B _u  \Big \}.$$
Thus, by recalling the definition of $S_ { B} (r),$ and taking $\widehat B$ in place of $B$ and $ \rho (r)$ in place of $r,$ we get:
\begin{equation}
\rho (\tau _r)= \inf \Big \{ s \ge \rho(r): \widehat B_s \ge M^{ \widehat B} _{\rho (r) } \Big \}= S_ {\widehat B} (\rho (r)) + \rho (r) ,
\end{equation}
where $M^{ \widehat B} _r = \max _{u \in [0,r]} \widehat B_u \ .$ Therefore, $\tau  _r = \rho ^{-1}
\left ( S _{ \widehat B} (\rho (r)) + \rho (r) \right ),$
and $ S(r) = \tau_r-r = \rho ^{-1} \left (S _{\widehat B} (\rho (r)) + \rho (r) \right ) -r .$ Finally:
$$ P \{ S(r) \le t \}= P \left \{ \rho ^{-1} (S _ {\widehat B} (\rho (r)) + \rho (r)) \le t+r \right \} $$
\begin{equation}
= P \left \{ S _ {\widehat B} (\rho (r)) + \rho (r) \le \rho (t+r) \right \} =
P \left \{ S _ {\widehat B} (\rho (r)) \le \rho (t+r) - \rho (r)\right \},
\end{equation}
from which \eqref{result} follows, by using \eqref{distrtauBM};
formula \eqref{denresult} is obtained by taking the derivative with respect to $t.$ \par\noindent
(ii) Suppose that  $\rho (t)$ is not deterministic, and the bounds \eqref{boundsrho} hold; set:
\begin{equation}
 \tau _ {r, \alpha } = \inf \{ t \ge r: \widehat B ( \rho (t)) \ge \max _ { u \in [0, \alpha (r)] } \widehat B _u \} , \
\tau _ {r, \beta } = \inf \{ t \ge r: \widehat B ( \rho (t)) \ge \max _ { u \in [0, \beta (r)] } \widehat B _u \} .
\end{equation}
As easily seen, one has:
\begin{equation}
\tau _ {r, \alpha } \le \tau_r \le \tau _ {r, \beta } \ ,
\end{equation}
that implies:
\begin{equation} \label{boundsfortaus}
\rho ( \tau _ {r, \alpha }) \le \rho (\tau_r) \le \rho ( \tau _ {r, \beta }).
\end{equation}
Moreover, since $\rho ( \tau _ {r, \alpha }) = \inf \{ s > \rho (r) : \widehat B _s \ge \max _ { u \in [0, \alpha (r)] } \widehat B _u \},$
and
$\rho ( \tau _ {r, \beta }) = \inf \{ s > \rho (r) : \widehat B _s \ge  \max _ { u \in [0, \beta (r)] } \widehat B _u \},$
we have :
\begin{equation}
\inf \{ s > \alpha (r) : \widehat B _s \ge \max _ { u \in [0, \alpha (r)] } \widehat B _u \} \le \rho ( \tau _ {r, \alpha })
\end{equation}
and
\begin{equation}
\rho ( \tau _ {r, \beta }) \le
\inf \{ s > \beta (r) : \widehat B _s \ge \max _ { u \in [0, \beta (r)] } \widehat B _u \} .
\end{equation}
Thus, recalling the definition of $S_{\widehat B} (r),$ from \eqref{boundsfortaus} we get:
\begin{equation}
S_{\widehat B} ( \alpha (r)) + \alpha (r)  \le \rho (\tau _r) \le
S_{\widehat B} ( \beta (r)) + \beta (r),
\end{equation}
and so:
\begin{equation}
 \rho ^ {-1} (  S_{\widehat B} ( \alpha (r)) + \alpha (r)) \le \tau_r \le \rho ^ {-1} (  S_{\widehat B} ( \beta (r)) + \beta (r)) ,
 \end{equation}
where $\rho ^ {-1}$ is the ``inverse'' of the random function $\rho,$ which is defined by $\rho ^ {-1} (s) := \inf \{t >0: \rho (t) >s  \}.$ Since \eqref{boundsrho} implies that $\beta ^ {-1} (s) \le \rho ^ {-1} (s) \le \alpha ^ {-1} (s) ,$ we obtain:
\begin{equation}
 \beta ^ {-1} (  S_{\widehat B} ( \alpha (r)) + \alpha (r)) \le \tau_r \le \alpha ^ {-1} (  S_{\widehat B} ( \beta (r)) + \beta (r)) .
 \end{equation}
Therefore:
\begin{equation} \label{boundtoSr}
\beta ^ {-1} (  S_{\widehat B} ( \alpha (r)) + \alpha (r)) -r \le S(r) \le \alpha ^ {-1} (  S_{\widehat B} ( \beta (r)) + \beta (r)) -r .
\end{equation}
From the first inequality in \eqref{boundtoSr}, it follows that:
\begin{equation}
 P \{ S(r) \le t \} \le P \{ S_{\widehat B} ( \alpha (r)) \le \beta (t+r) - \alpha (r) \} =
\frac 2 \pi \arctan \left ( \sqrt {\frac {\beta (t+r) - \alpha (r)} {\alpha (r)} } \ \right ) ,
\end{equation}
that proves \eqref{bound2}. Moreover, if
there exists $\bar t >0$ such that $\alpha( r+ \bar t ) \ge \beta (r),$ we have
$\alpha (r+ t ) \ge \beta (r) $ for $t > \bar t ,$ because $\alpha (t)$ is increasing; then, for $t > \bar t$ the second inequality in \eqref{boundtoSr} implies:
\begin{equation}
 P \{ S(r) \le t \} \ge P \{ S_{\widehat B} ( \beta (r)) \le \alpha (t+r) - \beta (r) \} =
\frac 2 \pi \arctan \left ( \sqrt {\frac {\alpha (t+r) - \beta (r)} {\alpha (r)} } \ \right ) ,
\end{equation}
which proves \eqref{bound1}.
The condition $\alpha (t +r) \ge \beta (r)$ is necessary so that $S_ {\widehat B} ( \beta (r)) \ge 0;$
of course, if  a value $\bar t>0,$ such that $\alpha (r+ \bar t) \ge \beta (r) ,$ does not exist, the inequality \eqref{bound1} loses meaning, because the square root is not defined.
\hfill $\Box$

\begin{Remark} {\rm
Notice that \eqref{result} is independent of the scale function $w;$ in particular, if $X$ is conjugated to BM via the function
$v,$ being $\rho (t)= t,$  one obtains that the distribution of $S (r)$ is the same as that of $S _B (r),$  given by \eqref{distrtauBM}.
\par\noindent
For an example of diffusion $X$ for which $\rho$ is not deterministic, but satisfies the bounds \eqref{boundsrho} with $\alpha (t)$ close to $\beta (t),$  see Example 4 of (Abundo, 2012).
}
\end{Remark}

\begin{Remark} {\rm
Let us suppose that $\rho (t)$ is deterministic and there exists $\gamma >0$ such that $\rho (t) \sim const \cdot t^ \gamma , $ as $t \rightarrow + \infty;$ then, from \eqref{denresult} it easily follows that
$E( S (r)) < + \infty ,$ provided that $\gamma >2 .$ This is not the case of BM, because $\rho (t) =t;$ instead, it holds e.g. for integrated BM
\ $X(t)= \int _0 ^t B_s ds,$ being $X(t)= \widehat B (\rho (t))$ with
$\rho (t) = t^3/3 $ (see (Abundo, 2015), (Abundo, 2013)).
}
\end{Remark}

\begin{Remark} {\rm
Let us suppose that $\rho (t)$ is deterministic; for given $r_1, r_2$ with $0\le r_1 < r_2,$ set $M_{[r_1, r_2]} := \max _{t \in [r_1, r_2] } X(t),$ and $S(r_1, r_2) = \inf \{ t \ge r_2 : X(t)=
M_{[r_1, r_2] } \} -r_2.$ Then, by using the arguments of the proof of Theorem \ref{teorema} and the Remark of (Papanicolaou, 2016) which refers to BM,
one obtains:
\begin{equation}
P \{S(r_1, r_2) \le t \}= \frac 2 \pi \arctan \left ( \sqrt {\frac {\rho (t+r_2) - \rho (r_2)} {\rho (r _2) - \rho (r_1)} } \ \right ), \ t \ge 0.
\end{equation}
}
\end{Remark}

\begin{Remark} {\rm Let us suppose that $\rho (t)$ is deterministic, and for $r>0,$ define $L(r) = \min _{t \in [0,r]} X(t) ,$ and
$U(r)= \inf \{t \ge r: X(t) \ge L(r) \} -r .$  Then, by using the fact, proved in (Papanicolaou, 2016), that $S_B(r)$ and $U_B(r)
:=  \inf \{t \ge r: B^ \eta _t \ge L_B(r) \} -r$ have the same distribution
(here $L_B(r) =  \min _{t \in [0,r]} B^ \eta _t ),$ and
by arguments analogous to those in the proof of Theorem \ref{teorema}, we conclude that also
$U(r)$ and $S(r)$ have a common distribution.

}
\end{Remark}

\section{Conclusions and final remarks}
We have considered a time-homogeneous one-dimensional diffusion process $X$ in an interval $I \subset \mathbb{R},$ driven by
the SDE \eqref{eqdiffu};
for a given time $r,$  under suitable conditions, we have found the distribution of the
time $S(r)$ required (after $r)$ for
$X(t)$  to exceed its maximum
$M_r = \max _{t \in [0,r]} X(t)$ on the interval $[0,r],$
generalizing the result in (Papanicolaou, 2016), which refers to
BM. Indeed, we have reduced $X$ to BM (see \eqref{representation}) by a space transformation,
given by the scale function $w(x),$ and a random time-change $\rho (t),$ under the assumption that
 $\rho ( + \infty) = + \infty $ ($\rho (t)$ is the quadratic variation of the space-transformed process).
Thus, we have shown that, when $\rho (t)$ is deterministic, $S(r)$ follows a compound arctangent law; note that, in this case
$X(t)$ solves also
the SDE:
\begin{equation}
dX(t)= - \frac {\rho ' (t) w'' (X(t)) } { 2 (w' (X(t)))^3 } dt + \frac {\sqrt { \rho '(t)} } { w' (X(t))}  d \widetilde B_t \ ,
\end{equation}
where  $w'(x)$ and $w''(x)$ denote first and second derivative of $w(x),$ and $\widetilde B$ is BM (see also (Abundo, 2017)).
The class of processes $X,$  for which the distribution of $S(r)$ has been found, includes,
besides diffusions conjugated to BM (see e.g. (Abundo, 2012)),
integrated BM and  integrated Ornstein-Uhlenbeck  process (see (Abundo, 2015), (Abundo, 2013)). \par
As a curiosity, we note that a number of results are known, which regard inverse trigonometric laws for some random times
associated to BM;
for instance, the density \eqref{densityofSB} appears as the conditional density of the second inter-passage time of BM
through a level, with the condition that the first-passage time is $r$ (see eq. (2.15) of (Abundo, 2016)). \par\noindent
The arc-sine law is valid for the time $\tau$ spent by BM on the positive half-line during the time interval
$[0,r],$ that is, $P( \tau \le t ) = \frac 2 \pi \arcsin \left ( \sqrt {\frac {t } { r}} \ \right ), \ t \in [0,r]$ (see (Levy, 1965)); moreover,
a compound arc-sine
law holds for the first
instant $\theta$ at which a diffusion $X$ of the form \eqref{representation}, with $\rho (t)$ deterministic, attains the maximum in the interval $[0,r],$ namely
$P( \theta \le t )= \frac 2 \pi \arcsin  \left (\sqrt { \frac {\rho (t) } { \rho (r)} } \ \right ), \ t \in [0,r]$ (see (Abundo, 2006),
 and (Levy, 1965) in the case of BM, i.e. $\rho (t)=t).$
\newpage

\begin{large}
\noindent {\bf References}
\end{large}

\bigskip

\noindent Abundo, M., 2017.
The mean of the running maximum of an
integrated Gauss-Markov process and the connection with its
first-passage time.
Stochastic Anal. Appl. 35:3, 499-510, DOI: 10.1080/07362994.2016.1273784
\par\noindent
Abundo, M., 2016. On the excursions of drifted Brownian motion and the successive passage times of Brownian motion.
Physica A 457, 176--182. \par\noindent
Abundo, M., 2015. On the first-passage time of an integrated
Gauss-Markov process. Scientiae Mathematicae Japonicae Online
e-2015, 28,  1--14.\par\noindent
Abundo, M., 2013. On the
representation of an integrated Gauss-Markov process. Scientiae
Mathematicae Japonicae Online e-2013, 719–-723. \par\noindent
Abundo, M., 2012. An inverse first-passage problem for
one-dimensional diffusions with random starting point. Statist.
Probab. Lett. 82, 7–-14. \par\noindent
Abundo, M., 2006. The arc-sine law for the first instant at which a diffusion process equals the ultimate value of a functional.
Int. J. Pure Appl. Math., 30 (1), 13--22.  \par\noindent
Dassios, A. and Lim, J.W., 2017.
Methodol Comput Appl Probab Online first, 25 Jan 2017,
doi:10.1007/s11009-017-9542-y
\par\noindent
Ikeda, N. and Watanabe, S., 1981.
Stochastic differential equations and
diffusion processes.
North-Holland Publishing Company.
\par\noindent
Karlin, S. and Taylor, H.M., 1975.
A second course in stochastic processes.
Academic Press, New York.
\par\noindent
Levy, P., 1965. Processus Stochastiques et Mouvement Brownien. Gauthier-Villars, Paris. \par\noindent
Papanicolaou, V.G., 2016.
An arctangent law.
Statist. Probab. Lett. 116, 62-–64.
\par\noindent
Revuz, D. and Yor, M., 1991.
Continous martingales and Brownian motion.
Springer-Verlag, Berlin Heidelberg.

\end{document}